\documentclass{ifacconf}

\usepackage{natbib}

\usepackage[T1]{fontenc}
\usepackage{enumerate}
\usepackage[latin9]{inputenc}
\usepackage{amsmath}
\usepackage{amssymb}
\usepackage{graphicx}
\usepackage{dsfont}

\makeatletter

\pdfpageheight\paperheight
\pdfpagewidth\paperwidth

\@ifundefined{showcaptionsetup}{}{%
 \PassOptionsToPackage{caption=false}{subfig}}
\usepackage{subfig}
\makeatother

\newcommand{\spa}[1]{\text{lin}\left(#1\right)}

\newcommand{\I}{\mathcal{I}}
\newcommand{\La}{\mathcal{L}}
\newcommand{\E}{\mathcal{E}}
\newcommand{\Prob}{\mathbb{P}}
\newcommand{\Ht}{\mathcal{H}}
\newcommand{\V}{\mathcal{V}}

\newcommand{\G}{\mathcal{G}}
\newcommand{\Ex}{\mathbb{E}}
\renewcommand{\e}[2]{e_{#1}^{#2}}
\newcommand{\de}[2]{d_{#1}^{#2}}
\newcommand{\ve}[3]{v_{#1}^{#2,#3}}
\def\QED{~\rule[-1pt]{5pt}{5pt}\par\medskip}
\newtheorem{defin}{Definition}[section]
\newtheorem{rema}{Remark}[section]
\newtheorem{lemma}{Lemma}[section]

\newtheorem{thrm}{Theorem}[section]
\newtheorem{problem}{Problem}[section]
\newcommand{\set}[1]{\left\{ #1 \right\}}

\begin{document}
\begin{frontmatter}

\thanks[footnoteinfo]{This research was in part supported by NSF Career award CCF-1350590, and by the NSF, AFOSR, ARPA-E, and the Institute for Collaborative Biotechnologies through grant W911NF-09-0001 from the U.S. Army Research Office. The content does not necessarily reflect the position or the policy of the Government, and no official endorsement should be inferred.}
\thanks[footnoteinfo2]{A. Lamperski was supported by a grant from the Whitaker International Scholars Program.}

\title{Optimal Two Player LQR State Feedback With Varying Delay\thanksref{footnoteinfo}\thanksref{footnoteinfo2}}

\author[First]{Nikolai Matni} 
\author[Second]{Andrew Lamperski} 
\author[First]{John C. Doyle}

\address[First]{Department of Control and Dynamical Systems, California Institute of Technology, Pasadena, CA, 91125 USA (e-mail: \{nmatni,doyle\}@caltech.edu)}                                              
\address[Second]{Department of Engineering, University of Cambridge, Cambridge CB2 1PZ, UK
(e-mail: a.lamperski@eng.cam.ac.uk)}

\begin{abstract}
This paper presents an explicit solution to a two player distributed
LQR problem in which communication between controllers occurs across
a communication link with varying delay. We extend known dynamic
programming methods to accommodate this varying delay, and show
that under suitable assumptions, the optimal control actions are linear in their information, and that the resulting controller has piecewise linear 
dynamics dictated by the current effective delay regime. 
\end{abstract}

\end{frontmatter}

\section{Introduction}
\label{sec:intro}
In the past decade, optimal decentralized controller synthesis has seen an explosion of advances at
the theoretical, algorithmic and practical levels. We provide a brief
survey of the more directly relevant results to our paper in the following,
and refer the reader to the tutorial paper by \cite{MMRY12}
for a timely presentation of the current state of the
art in optimal decentralized control subject to information constraints. 

A particular class of decentralized control problems that has received a significant amount of attention is that of optimal $\mathcal{H}_{2}$ (or LQG) control subject to delay constraints. In this case, the information constraints can be interpreted as arising from a communication graph, in which edge weights between nodes correspond to the delay required to transmit information between them.  For the special case of the one-step delay information sharing pattern, the $\mathcal{H}_{2}$ problem was solved in the 1970s using dynamic programming (\cite{SA74,KS74,Y75}). For more complex delay patterns, sufficient statistics are not easily identified, making extensions beyond the state feedback case (\cite{LD11,LD12}) difficult, although semi-definite programming (SDP) (\cite{R06,G06}), vectorization (\cite{RL06}), and spectral factorization (\cite{LDXX}) based solutions do exist.  It is worth noting that for specific systems, sufficient statistics and a generalized separation principle have been identified and successfully applied, as in the work by \cite{FGJ12}.   Furthermore, recent work by \cite{NMT11,NMT13} provides dynamic programming decompositions for the general delayed sharing model. 

An underlying assumption in all of the above is that information,
albeit delayed, can be transmitted \emph{perfectly} across a communication
network with a \emph{fixed} delay. A realistic communication network,
however, is subject to data rate limits, quantization, noise and packet
drops \textendash{} all of these issues result in possibly varying
delays (due to variable decoding times) and imperfect transmission
(due to data rate limits/quantization).  The assumption that these delays are fixed necessarily
introduces a significant level of conservatism in the control design procedure.
In particular, to ensure that the delays under which controllers 
exchange information do not vary, worst case delay times must be used
for control design, sacrificing performance and robustness in the process.

These issues have been addressed by the networked control
systems (NCS) community, leading to a plethora of results for channel-in-the
loop type problems: see the recent survey  by \cite{Hesp}, and the references therein.  Some of the more relevant results from this field include the work by \cite{GSHM08} and  \cite{GSC10}, which address optimal
LQG control of a single plant over a packet dropping channel.  Very few results exist, however, that seek to combine NCS and decentralized optimal control.  A notable exception is the work by \cite{CL11}, in which an explicit state space solution to a sparsity constrained two-player decentralized LQG state-feedback problem over a TCP erasure channel is solved. 

We take a different view from these results, and suppress the underlying details of the communication network, and instead assume that packet drops, noise, and congestion manifest themselves to the controllers as varying delays.  In particular, we seek to extend the distributed state-feedback results of \cite{LD11,LD12} and \cite{LL12}
to accommodate varying delays. In addition to allowing for communication
channels to be more explicitly accounted for in the control design
procedure, the ability to accommodate varying delays provides
flexibility in the coding design aspect of this problem.

In this paper, we focus on a two plant system in which communication
between controllers occurs across a communication link with varying
delay. In \cite{MD_CDC13_jitter}, we solved a special case of this problem by extending the methods used in \cite{LD11} and \cite{LD12}.  Here, we use a variant of the dynamic programming methods in \cite{LL12} to
accommodate this varying delay, and show that under suitable assumptions, the optimal
control actions are linear in their information, and that the resulting controller has piecewise linear 
dynamics dictated by the current effective delay regime. 

This paper is structured as follows: in Section \ref{sec:problem} we fix notation,
and present the problem to be solved in the paper. Section \ref{sec:main} introduces the concepts of effective delay, partial nestedness (c.f. \cite{HC72}) and a system's information graph (c.f. \cite{LL12}) before presenting our main result.  Section \ref{sec:derivation} derives the optimal control actions and controller, and Section \ref{sec:conclusion} ends with conclusions and directions for future work.  Proofs of all intermediary results can be found in the Appendix.

\section{Problem Formulation}
\label{sec:problem}
\subsection{Notation}
For a matrix partitioned into blocks \[ M=\left[\begin{array}{ccc} M_{11} & \cdots & M_{1N}\\ \vdots & \ddots & \vdots\\ M_{N1} & \cdots & M_{NN} \end{array}\right] \] and $s,v\subset\{1,\dots,N\},$ we let $M^{s,v}=(M_{ij})_{i\in s,j\in v}$. 

For example \[ M^{\{1,2,3\}\{1,2\}}=\left[\begin{array}{cc} M_{11} & M_{12}\\ M_{21} & M_{22}\\ M_{31} & M_{32} \end{array}\right]. \] We denote the sequence $x_{t_{0}},...,x_{t_{0}+t}$ by $x_{t_{0}:t_{0}+t}$, and given the history of a random process $r_{0:t}$, we denote the conditional probability of an event $\mathcal{A}$ occurring given this history by $\mathbb{P}_{r_{0:t}}(\mathcal{A})$.  If $\mathcal{Y}=\{y^1,\dots,y^M\}$ is a set of random vectors (possibly of different sizes), we say that $z \in \spa{\mathcal{Y}}$ if there exist appropriately sized real matrices $C^1,\dots,C^M$ such that $z=\sum_{i=1}^M C^iy^i$.


\subsection{The two-player problem}

This paper focuses on a two plant system with physical propagation
delay of $D$ between plants, and stochastically varying communication
delays $d^{i}_t\in\{0,\dots,D\}$ -- to ease notation, we let $d_t:=(d^{1}_t,d^2_t)$. We impose some additional assumptions on the stochastic process $d_t$ in Section \ref{sec:main} such that the infinite horizon solution is well defined.

The dynamics of the sub-system $i$ are then captured by the following difference equation:
\begin{equation}
\begin{array}{ccc}
x^{i}_{t+1} & = & A_{ii}x^{i}_t+A_{ij}x^{j}_{t-(D-1)}+B_{i}u^{i}_t+w^i_t
\end{array}
\label{eq:local}
\end{equation}
with mutually independent Gaussian initial conditions and noise vectors
\begin{equation}
x^i_0 \sim \mathcal{N}(\mu_0^i,\Sigma^i_0),\ \ \ \ w^i_t \sim \mathcal{N}(0,W^i_t)
\label{eq:noise_stats}
\end{equation}  

We may describe the information available to controller $i$ at time $t$, denoted by $\I^i_t$, via the following recursion:
\begin{equation}
\begin{array}{rcl}
\I^i_0&=&\{x^i_0\}  \hfill{}\\
\I^i_{t+1} &=& \I^i_t\cup\{x^i_{t+1}\}\cup\{x_k^j\,:\, 1 \leq k \leq t+1-d^j_{t+1}\} \hfill{}
\end{array}
\label{eq:info_recursion}
\end{equation}

The inputs are then constrained to be of the form
\begin{equation}
\begin{array}{rcl}
u^{i}_t & = & \gamma^i_t(\I^i_t)
\label{eq:input_cons}
\end{array}
\end{equation}
for Borel measurable $\gamma^i_t$.

\begin{figure}[t]
\begin{center}
\includegraphics[width=2.5in]{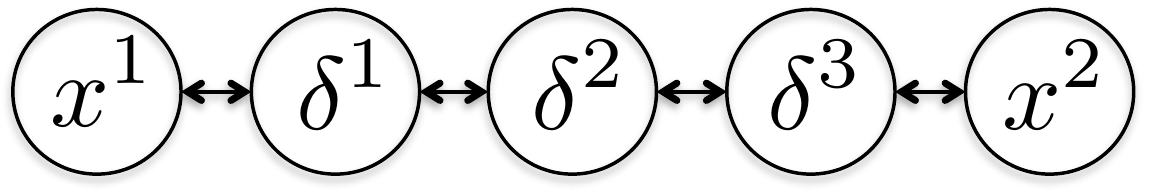}
\caption{The distributed plant considered in (\ref{eq:statespace}), shown here for $D=4$.  Dummy nodes $\delta^i_t$, $i=1,\dots,D-1$, as defined by (\ref{eq:deltas}), are introduced to make explicit the propagation delay of $D$ between plants.   \label{fig:plant}}
\end{center}
\end{figure}

In order to build on the results in \cite{LL12}, we model the two
plant system as a $D+1$ node graph, with ``dummy delay''
nodes introduced to explicitly enforce the propagation delay between
plants. Specifically, letting 
\begin{equation}
\delta^i_t=\left[\begin{array}{c} x^{1}_{t-i} \\ x^{2}_{t-(D-i)}\end{array}\right], \ i=1,\dots,D-1
\label{eq:deltas}
\end{equation}
where $\delta^i$ is the state of the $i^\text{th}$ dummy node, we obtain the following state space representation for the system
\begin{equation}
x_{t+1}=Ax_t + Bu_t +  w_t
\label{eq:statespace}
\end{equation}
where, to condense notation, we let
\begin{equation}
\begin{array}{lcr}
x = \begin{bmatrix} x^1 \\ \delta^1 \\ \vdots \\ \delta^{D-1} \\ x^2 \end{bmatrix} &
u = \begin{bmatrix} u^1 \\ 0 \\ \vdots \\ 0 \\ u^2 \end{bmatrix} &
w = \begin{bmatrix} w^1 \\ 0 \\ \vdots \\0 \\ w^2 \end{bmatrix},
\end{array}
\end{equation}
 and $A$ and $B$ are such that (\ref{eq:statespace}) is consistent with (\ref{eq:local}) and (\ref{eq:deltas}).  The physical topology of the plant is illustrated in Figure \ref{fig:plant}.

\begin{figure*}[t]
\includegraphics[width=\textwidth]{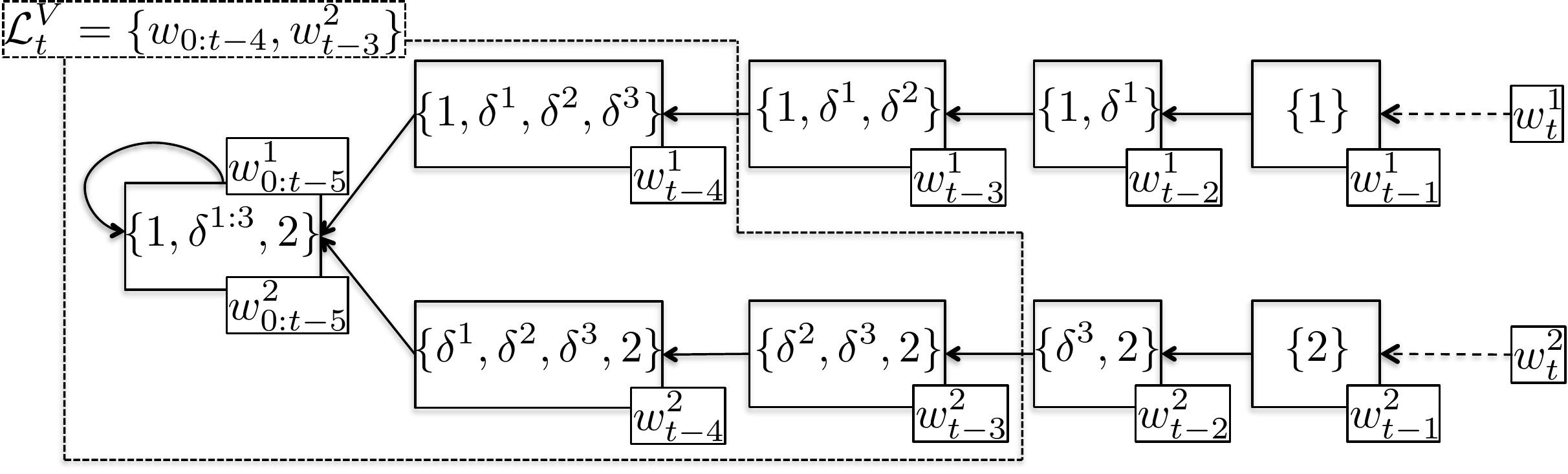}
\caption{The information graph $\G=(\V,\E)$, and label sets $\{\La^s_t\}_{s\in\V}$, for system \eqref{eq:statespace}, shown here for $D=4$, and $e_t=(3,2)$. Notice that: (i) for each $(r,s)\in \E$, with $|r|<D+1$, we have that $|s|=|r|+1$, (ii) that $|s|$ corresponds exactly to how delayed the information in the label set is, and (iii) that $\La^V_t$ contains all of the information at nodes s.t. $|s|> e^i_t,\,s\ni i$.  We also see that the graph is naturally divided into two branches, with each branch corresponding to information pertaining to a specific plant.}
\label{fig:info_graph}
\end{figure*}

\begin{problem}
Given the linear time invariant (LTI) system described by (\ref{eq:local}), (\ref{eq:deltas}) and (\ref{eq:statespace}), with disturbance statistics (\ref{eq:noise_stats}), minimize the infinite horizon expected cost
\begin{equation}
 \lim_{N\rightarrow\infty}\frac{1}{N}\mathbb{E}\left[\sum_{t=1}^{N}x(t)^{T}Qx(t)+u(t)^{T}Ru(t)\right] 
\label{eq:infcost}
\end{equation}
subject to the input constraints (\ref{eq:input_cons}). 
\label{problem}
\end{problem}

The weight matrices are assumed to be partitioned into blocks of appropriate dimension, i.e. $Q=(Q_{ij})$ and $R=(R_{ij})$, conforming to the partitions of $x$ and $u$. We assume $Q$ to be positive semi-definite, and $R$ to be positive definite, and  in order to guarantee existence of the stabilizing solution to the corresponding Riccati equation, we assume $(A,B)$ to be stabilizable and $(Q^{\frac{1}{2}},A)$ to be detectable.

\section{Main Result}
\label{sec:main}

\subsection{Effective delay}


The information constraint sets \eqref{eq:info_recursion} are defined in such a way that controllers do not forget information that they have already received.  This leads to the $x^j$ component of the information set $\I^i_t$ being a function of the \emph{effective} delay seen by the controller, as opposed to the current delay value of the communication channel $\de{t}{j}$.

\begin{defin}
Let 
\begin{multline}
\e{t}{j}:=\min \{d^j_{t},d^j_{t-1}+1,d^j_{t-2}+2,\dots,\\d^j_{t-(D-2)}+(D-2),d^j_{t-(D-1)}+(D-1)\}
\label{eq:eff_del}
\end{multline}
be the \emph{effective delay} in transmitting information from controller $j$ to controller $i$.
\end{defin}

\begin{lemma}
The information set available to controller $i$ at time $t$ may be written as
\begin{equation}
\I^i_t = \I^i_{t-1} \cup \{x^i_t\} \cup \I^j_{t-\e{t}{j}}
\label{eq:alt_info}
\end{equation}
\label{lem:alt_info}
\end{lemma}
\begin{pf}
See Appendix.
\end{pf}

In order to ensure that the infinite horizon solution is well defined, we assume that that the stochastic delay process $d_t$ induces an effective delay process such that
\begin{equation}
\lim_{T\to\infty}\frac{1}{T} \sum_{t=0}^{T-1}\Prob_{d_{0:t}}\left(e^i_{t+1}\leq d\right)
\label{eq:stationarity}
\end{equation} 
exists for any integer $d$. 
%


\subsection{Partial Nestedness}
Here we show that the information constraints (\ref{eq:input_cons}) and system (\ref{eq:statespace}) are \emph{partially nested} (c.f. \cite{HC72}), and hence that the optimal control policies $\gamma^i_t$ are linear in their information set.  

\begin{defin}
A system (\ref{eq:statespace}) and information structure (\ref{eq:input_cons}) is partially nested if, for every admissible policy $\gamma$, whenever $u^i_\tau$ affects $\I^j_t$, then $\I^i_\tau \subset \I^j_t$.
\end{defin}
\begin{lemma} (see \cite{HC72})
Given a partially nested information structure, the optimal control law that minimizes a quadratic cost of the form (\ref{eq:infcost}) exists, is unique, and is linear.
\label{lem:linctrl}
\end{lemma}

Using partial nestedness, the following lemma shows that the optimal state and input lie in the linear span of $\I^i_t$ and $\Ht_t$, where $\Ht_t$ is the \emph{noise history} of the system given by
\begin{equation}
\Ht_t = \{x_0,w_{0:t-1}\}
\label{eq:noisehist}
\end{equation}
\begin{lemma}
The system (\ref{eq:statespace}) and information structure (\ref{eq:input_cons}) is partially nested, and for any linear controller, we have that
\begin{equation}
x^i_t,u^i_t \in \spa{\I^i_t},\ \ \ \ \ x_t,u_t\in \spa{\Ht_t}
\end{equation}
\label{lem:partial}
\end{lemma}
\begin{pf}
See Appendix.
\end{pf}

\subsection{Information Graph and Controller Coordinates}

Lemma \ref{lem:partial} indicates that each $\I^i_t$ is a  subspace of $\Ht_t$: in this section, we exploit this observation to define pairwise independent controller coordinates.  An explicit characterization of these subspaces is given in Section \ref{sec:derivation}.  

We begin by defining the \emph{information graph}, as in \cite{LL12}, associated with system (\ref{eq:statespace}) by $\G = (\V,\E)$, with
\begin{equation}
\begin{array}{rcl}
\V &:=& \set{ \set{1},\set{1,\delta^1},\dots,\set{1,\delta^1,\dots,\delta^{D-1}} } \cup \\
&& \set{ \set{2},\set{\delta^{D-1},2},\dots,\set{\delta^1,\dots,\delta^{D-1},2} } \cup V\\
\E &:=& \set{ (r,s) \in \V \times \V \, : \, |s| = |r|+1}\cup\set{(V,V)}
\end{array}
\label{eq:graph}
\end{equation}
where $V:=\set{1,\delta^1,\dots,\delta^{D-1},2} $. For the case of $D=4$, the graph $\G$ is
  illustrated in Figure \ref{fig:info_graph}. 
  
Before proceeding, we define the following sets, which will help us state the main result.  Let
\begin{equation}
\begin{array}{rcl}
\ve{t}{i}{+} &:=&\{s\in\V\backslash{}V\,|\,i\in s, \, |s|\geq \e{t}{i} \} \\
\ve{t}{i}{++}  &:=&\{s\in\V\backslash{}V\,|\,i\in s, \, |s| > \e{t}{i} \}
\end{array}
\label{eq:helpers}
\end{equation}
and similarly define $\ve{t}{i}{-} $ and $\ve{t}{i}{--} $ as in \eqref{eq:helpers}, but with the (strict) inequality reversed.

\begin{thrm}
Consider Problem \ref{problem}, and let $\G(\V,\E)$ be the associated information graph.  Let
\begin{equation}
\begin{array}{rcl}
X^V &=& Q + A^\top X^V A + A^\top X^VBK^V \\
K^V &:=& -\left(R+B^\top X^VB\right)^{-1}B^\top A,
\end{array}
\end{equation}
be the stabilizing solution to the discrete algebraic Riccati equation, and the centralized LQR gain, respectively.  Now, assume that $X^s$ is given, and let $r\neq s \in\V$ be the unique node such that $(r,s)\in\E$. Define the matrices
\begin{equation}
\begin{array}{ll}
\Lambda^r = Q^{rr} +  p^r (A^{Vr})^\top X^VA^{Vr} + q^r (A^{sr})^\top X^sA^{sr} \hfill{}\\
\Psi^r = R^{rr} + p^r (B^{Vr})^\top X^VB^{Vr} + q^r (B^{sr})^\top X^sB^{sr} \hfill{} \\
\Omega^r = p^r (A^{Vr})^\top X^VB^{Vr} + q^r (A^{sr})^\top X^sB^{sr} \hfill{} \\
X^r = \Lambda^r + \Omega^rK^r \hfill{} \\
K^r = -\left(\Psi^r\right)^{-1}\left(\Omega^r\right)^\top
\end{array}
\end{equation}
where $p^r$ is given by
\begin{equation}
p^r := \lim_{T\to\infty}\frac{1}{T} \sum_{t=0}^{T-1}\Prob_{d_{0:t}}\left(r\in \ve{t+1}{i}{++}\right)
\label{eq:probs}
\end{equation}
and $q^r = 1- p^r$.

The optimal control decisions then satisfy
\begin{equation}
\begin{array}{rcl}
\zeta^V_{t+1} &=& A\zeta^V_t + B\varphi^V_t +\\ && \indent \sum_{i=1}^2 \sum_{r\in\ve{t+1}{i}{++}}(A^{Vr}\zeta^r_t + B^{Vr}\varphi^r_t)\\
\zeta^s_{t+1} &=& \begin{cases} A^{sr}\zeta^r_t + B^{sr}\varphi^r_t \ \text{if } s\in \cup_i\ve{t+1}{i}{-}, \,(r,s)\in\E \\ 
 0 \ \ \ \ \ \text{otherwise}
\end{cases}\\
\zeta^i_{t+1} &=& w^i_t \\
\zeta^i_{0} &=& x^i_0 \\
u^i_t &=&  \varphi^V_t + \sum_{s\in\ve{t}{i}{-}} I^{V,s}\varphi^s_t \\
\varphi^r_t & = & K^r\zeta^r_t
\end{array}
\label{eq:dynamics_main} 
\end{equation}
and the corresponding infinite horizon expected cost is
\begin{equation}
\sum_{i=1}^2\text{Trace}\left(X^{\{i\}}W^i\right)
\end{equation}
\end{thrm}
\begin{pf}
See Section \ref{sec:derivation}.
\end{pf}
\begin{rema}
Notice that the global action taken based on $\zeta^V$ must be taken simultaneously by \emph{both} players.  In other words, it is assumed that an acknowledgment mechanism is in place such that $e_t$ is known to both players; relaxing this assumption will be the subject of future work.
\end{rema}
\begin{rema}
The probabilities $p^r$ and $q^r$ can be computed directly if we assume the $\{d_t\}$ to be independently and identically distributed.  In this case, $e^j_t$  evolves according to an irreducible and aperiodic Markov chain with transition probability matrix computable directly from the definition of effective delay and the pmf of $d_t$.  As such, $p^r$ and $q^r$ can be computed from the chain's stationary distribution, which is guaranteed to exist.  Future work will explore what additional distributions on $d_t$ will lead to closed form expressions for $p^r$ and $q^r$.  Failing the existence of closed form expressions for these asymptotic distributions, computing estimates via simulation should be a feasible option for many interesting delay processes.
\end{rema}

\section{Controller Derivation}
\label{sec:derivation}
\subsection{Controller States and Decoupled Dynamics}

As mentioned previously, each $\I^i_t$ is a  subspace of $\Ht_t$: in this section, we aim to explicitly characterize these subspaces by assigning label sets $\{\La^s_{0:t}\}_{s\in\V}$ to the graph $\G=(\V,\E)$ as defined by (\ref{eq:graph}).  In particular, they are defined recursively as:
\begin{equation}
\begin{array}{rcl}
\La^s_0 &=&\emptyset, \ \text{for } |s|>1\\
\La^{i}_0 &=& \{x^i_0\} \\
\La^{i}_{t+1} &=& \{w^i_{t}\} \\
\La^s_{t+1} &=& \La^{r}_t,\ \text{for }(r,s)\in \E,\ 1<|s|<D+1 \\
\La^V_{t+1} &=& \La^V_t\displaystyle\cup_i\cup_{s \in \ve{t+1}{i}{+}}\La^s_t
\end{array}
\label{eq:labelrecursions}
\end{equation}
where we have let $\cup_i$ denote $\cup_{i=1}^2$ to lighten notational burden.  An example of these label sets for the case of $D=4$ is illustrated in Figure \ref{fig:info_graph}.

Before delving in to the technical justification for these label sets, we provide some intuition.  The information graph $\G$ characterizes how the effect of noise terms spread through the system, and labels are introduced as a means of explicitly tracking this spreading.  As can be seen in Figure \ref{fig:info_graph}, for each $(r,s)\in \E$, with $|r|<D+1$, we have that $|s|=|r|+1$, and additionally, that $|s|$ measures exactly how delayed the information in the label set is.  We also see that the graph is naturally divided into two \emph{disjoint} branches, with each branch corresponding to information about a specific plant.  Finally, the label corresponding to the root node $V$ can be interpreted as the information available to both controllers -- this is reflected by its explicit dependence on the effective delay $\e{t}{i}$.

\begin{rema}
Note that in contrast to \cite{LL12}, the label sets as defined will in general \emph{not} be disjoint.  However, as will be made explicit in Lemma \ref{lem:partition}, an effective delay dependent subset of the label sets will indeed form a partition (i.e. a pairwise disjoint cover) of the noise history.
\end{rema}

We may now characterize the subspaces of $\Ht_t$ that are associated with each $\I^i_t$.  This characterization will be shown to depend on the effective delay $\e{t}{j}$ seen at node $i$, and will lead to an intuitive partitioning of both the state and the control input. 

We begin by pointing out the following useful facts that will be used repeatedly in the derivation to come
\begin{lemma}
Let $\ve{t}{i}{*}$, $*\in\{-,--\}$, be given as in (\ref{eq:helpers}).  Then, for a fixed $i$, we have that
\begin{equation}
\cup_{s\in \ve{t+1}{i}{-}}\La^s_{t+1} = \cup_{r \in \ve{t+1}{i}{--}} \La^r_t\cup\La^i_{t+1},
\end{equation}
and for integers $a,b\in\{0,\dots,D-1\}$
\begin{equation}
\cup_{ a<|s|\leq b+1}\La^s_{t+1} = \cup_{a \leq |r| \leq b} \La^r_{t}
\end{equation}
\label{lem:shrink}
\end{lemma}
\begin{pf}
Follows immediately by applying the recursion rules \eqref{eq:labelrecursions} and the fact that for each $(r,s)\in \E$, with $|r|<D+1$, we have that $|s|=|r|+1$.
\end{pf}

 \begin{lemma}
\label{lem:partition}
Consider the information graph $\G$ as defined in equation (\ref{eq:graph}), and the label sets defined as in \eqref{eq:labelrecursions}.
We then have that
\begin{enumerate}[(i)]
	\item For all $t\geq 0$, a subset of the labels form a partition of the noise history.  In particular, we have that
\begin{equation}
\Ht_t = \La^V_t \cup_i\cup_{s\in \ve{t}{i}{-}}\La^s_t
\label{eq:part}
\end{equation}
where the union is disjoint,  i.e. $\La^V_t\cap\La^s_t=\emptyset$ if $s\in\ve{t}{i}{-}$, and $\La^s_t\cap\La^{s'}_t = \emptyset$ for any $s \neq s', \, s,s' \in \cup_i \ve{t}{i}{-} $.
\item For $i=1,2$
\begin{equation}
\spa{\I^i_{t}} = \spa{\La^V_t \cup_{s\in \ve{t}{i}{-}}\La^s_t}.
\end{equation}
\end{enumerate}
\end{lemma}
\begin{pf}
See Appendix.
\end{pf}

\begin{rema}
Although the proof of this result is notationally cumbersome, it is mainly an exercise in bookkeeping.  The idea is illustrated in Figure \ref{fig:info_graph}: labels for nodes $v \neq V$ track the propagation of a disturbance through the plant, whereas the label for $V$ selects those labels corresponding to globally available information, as dictated by the effective delay.
\end{rema}

With the previous lemmas at our disposal, we may now write
\begin{equation}
\begin{array}{rcl}
x_t &=& \zeta^V_t + \sum_{i=1}^2\sum_{s\in\ve{t}{i}{-}} I^{V,s}\zeta^s_t \\
u_t &=& \varphi^V_t + \sum_{i=1}^2\sum_{s\in\ve{t}{i}{-}} I^{V,s}\varphi^s_t \\
\end{array}
\label{eq:decomp}
\end{equation}
where each $\zeta^s_t$, $\varphi^s_t \in \spa{\La^s_t}$.

We may accordingly derive update dynamics for these state and control components.

\begin{lemma}
If the control components are such that $\varphi^t_s \in \spa{\La^s_t}$, then the state components $\{\zeta^s_t\}$ satisfy the following update dynamics
\begin{equation}
\begin{array}{rcl}
\zeta^V_{t+1} &=& A\zeta^V_t + B\varphi^V_t +\\ && \indent \sum_{i=1}^2 \sum_{r\in\ve{t+1}{i}{++}}(A^{Vr}\zeta^r_t + B^{Vr}\varphi^r_t)\\
\zeta^s_{t+1} &=& \begin{cases} A^{sr}\zeta^r_t + B^{sr}\varphi^r_t \ \text{if } s\in \cup_i\ve{t+1}{i}{-}, \,(r,s)\in\E \\ 
 0 \ \ \ \ \ \text{otherwise}
\end{cases}\\
\zeta^i_{t+1} &=& w^i_t \\
\zeta^i_{0} &=& x^i_0
\end{array}
\label{eq:dynamics}
\end{equation}
\label{lem:dynamics}
\end{lemma}
\begin{pf}
See Appendix.
\end{pf}

In particular, notice that the dynamics \eqref{eq:dynamics} imply $\zeta^s_t = 0$ for all $s\in \cup_i\ve{t}{i}{++}$, allowing us to rewrite the decomposition for $x_t$ as
\begin{equation}
x_t = \sum_{s\in\V} I^{Vs}\zeta^s_t,
\label{eq:decomp2}
\end{equation}
where have simply added the zero valued state components to the expression in \eqref{eq:decomp}.

We now have all of the elements required to solve for the optimal control law via dynamic programming.

\subsection{Finite Horizon Dynamic Programming Solution}

Let $\gamma_t = \{\gamma_t^s\}_{s\in\V}$ be the set of policies at time $t$.  By Lemma \ref{lem:partial}, we may assume the $\gamma^s_t$ to be linear.  Define the cost-to-go
\begin{multline}
V_t(\gamma_{0:t-1}) = \hfill{}\\
\min_{\gamma_{t:T-1}} \Ex^{\gamma\times d} \left( \sum_{k=t}^{T-1} x_k^\top Q x_k + u_k^\top R u_k + x_T^\top Q_T x_T \right)
\end{multline}
where the expectation is taken with respect to the joint probability measure on $(x_{t:T},u_{t:T-1})\times(d_{t:T-1})$ induced by the choice of $\gamma = \gamma_{0:T-1}$ (note that the $d_t$ component is assumed to be independent of the policy choice).  

\begin{rema}
Following \cite{NMT13}, we adopt the common information formalism and define our cost-to-go function in terms of the control policy $\gamma$ to be chosen by a ``centralized coordinator.''  Noting that these policies can in fact be computed off-line and in a centralized manner (it is only their \emph{implementation} that requires measurement of the state components $\{\zeta^s\}$), this in effect reduces the dynamic programming argument to a standard full-information setting.
\end{rema}

Via the dynamic programming principle, we may iterate the minimizations and write a recursive formulation for the cost-to-go:
\begin{multline}
V_t(\gamma_{0:t-1}) = \hfill{}\\
\min_{\gamma_{t:T-1}} \Ex^{\gamma\times d}\left(  x_t^\top Q x_t + u_t^\top R u_t + V_{t+1}(\gamma_{0:t-1},\gamma_t) \right).
\label{eq:Vrec}
\end{multline}

We begin with the terminal time-step, $T$, and use the decomposition (\ref{eq:decomp2}) to obtain
\begin{equation}
V_T(\gamma_{0:T-1}) = \Ex^{\gamma\times d}\left( x_T^\top Q_T x_T\right) = \Ex^\gamma \sum_{s\in \V} (\zeta^s_T)^\top Q^{ss}_T (\zeta^s_T),
\end{equation}
where in the last step we have used the pairwise independence of the coordinates $\zeta^s_T$.  By induction, we shall show that the value function, for some $t\geq0$, always takes the form
\begin{equation}
V_{t+1}(\gamma_{0:t})=\Ex^{\gamma}\sum_{s\in\V}((\zeta^s_{t+1})^\top X^{s}_{t+1} (\zeta^s_{t+1})+c_{t+1}
\end{equation}
where $\{X^s_{t+1}\}_{s\in\V}$ is a set of matrices and $c_{t+1}$ is a scalar.  We now solve for $V_t(\gamma_{0:t-1})$ via the recursion \eqref{eq:Vrec}.  Given $e_t$, apply (\ref{eq:decomp2}) and the independence result to write
\begin{multline}
V_t(\gamma_{0:t-1}) = \\
\min_{\gamma_t} \Ex^{\gamma\times d} \left(\sum_{s\in\V} (\zeta^s_t)^\top Q^{ss}(\zeta^s_t) + (\varphi^s_t)^\top R^{ss}(\varphi^s_t) + \right. \\ 
\left. \sum_{s\in\V} (\zeta^s_{t+1})^\top X^{s}_{t+1}(\zeta^s_{t+1})+c_{t+1}\right)
\end{multline}

We now substitute the update equations \eqref{eq:dynamics}, average over $d_{t+1}$ and use independence to obtain
\begin{equation}
V_t(\gamma_{0:t-1}) = \min_{\gamma_t}\Ex^\gamma\left(\sum_{r \in \V} \begin{bmatrix} \zeta^r_t \\ \varphi^r_t \end{bmatrix}^\top \Gamma^r_t \begin{bmatrix} \zeta^r_t \\ \varphi^r_t \end{bmatrix}  + c_t\right)
\label{eq:solvefor}
\end{equation}
where $\Gamma^r_{0:T-1}$ and $c_{0:T-1}$ are given by:
\begin{multline}
\Gamma^r_t = \begin{bmatrix} Q^{rr} & 0 \\ 0 & R^{rr} \end{bmatrix} + \hfill{}\\ 
\Prob_{d_{0:t}}(r\in \ve{t+1}{i}{++})\begin{bmatrix} A^{Vr} & B^{Vr} \end{bmatrix}^\top X^V_{t+1}\begin{bmatrix} A^{Vr} & B^{Vr} \end{bmatrix} + \\
\Prob_{d_{0:t}}(r\in \ve{t+1}{i}{-})\begin{bmatrix} A^{sr} & B^{sr} \end{bmatrix}^\top X^s_{t+1}\begin{bmatrix} A^{sr} & B^{sr} \end{bmatrix}
\end{multline}
\begin{multline}
c_t = c_{t+1} + \sum_{i=1}^2\mathrm{Trace}\left(X^{\{i\}}_{t+1}W^i\right).\hfill{}
\label{eq:ct}
\end{multline}
The terminal conditions are  $c_T = 0$ and $\Gamma^r=Q^{rr}_T$, and $s$ is the unique node such that $(r,s)\in\E$.

Let $p^r_t:=\Prob_{d_{0:t}}(r\in \ve{t+1}{i}{++})$ and $q^r_t:=\Prob_{d_{0:t}}(r\in \ve{t+1}{i}{-})$, and introduce the following matrices:
\begin{equation}
\begin{array}{ll}
\Lambda^r_{t+1} = Q^{rr} +  p^r_t (A^{Vr})^\top X^V_{t+1}A^{Vr} + q^r_t (A^{sr})^\top X^s_{t+1}A^{sr} \hfill{}\\
\Psi^r_{t+1} = R^{rr} + p^r_t (B^{Vr})^\top X^V_{t+1}B^{Vr} + q^r_t (B^{sr})^\top X^s_{t+1}B^{sr} \hfill{} \\
\Omega^r_{t+1} = p^r_t (A^{Vr})^\top X^V_{t+1}B^{Vr} + q^r_t (A^{sr})^\top X^s_{t+1}B^{sr} \hfill{}
\end{array}
\end{equation}

Then each expression of the sum in \eqref{eq:solvefor} can be written as
\begin{equation}
(\zeta^r_{t})^\top\Lambda^r_{t+1}(\zeta^r_t) + (\varphi^r_t)^\top \Psi^r_{t+1}(\varphi^r_t) + 2(\zeta^r_t)^\top\Omega^r_{t+1}(\varphi^r_t).
\label{eq:expression}
\end{equation}

Due to the definitions of $\{\zeta^r\}$ and $\{\varphi^r\}$, it is clear that the terms \eqref{eq:expression} are pairwise independent and hence can be optimized independently.  Removing the information constraints, and optimizing over $\varphi^r_t$, we see that the optimal action is given by
\begin{equation}
\varphi^r_t = - \left(\Psi^r_{t+1}\right)^{-1}\left(\Omega^r_{t+1}\right)^\top\zeta^r_t
\label{eq:optimal_action}
\end{equation}
which, by construction, satisfies the information constraints $\I^i_t$.  Substituting this solution back in to (\ref{eq:expression}), we see that the matrices $X^r_t$ must satisfy
\begin{equation}
\begin{array}{rcl}
X^r_t &=& \Lambda^r_{t+1} + \Omega^r_{t+1}K^r_t \\
K^r_t &:=& - \left(\Psi^r_{t+1}\right)^{-1}\left(\Omega^r_{t+1}\right)^\top
\end{array}
\label{eq:recursion_and_gain}
\end{equation}

The finite horizon optimal cost is then given by
\begin{equation}
\begin{array}{rcl}
V_0 &=& \Ex \sum_{i=1}^2 (x^i_0)^\top X^{\{i\}}(x^i_0) + c_0\\
&=&\Ex \sum_{i=1}^2 (\mu_0^i)^\top X_0^{\{i\}}(\mu_0^i) + \text{Trace}\left( X_0^{\{i\}}\Sigma^i_0 \right) + c_0
\end{array}
\end{equation}
where $c_0$ can be computed according to \eqref{eq:ct} beginning with terminal conditions $c_T = 0$.

\subsection{Infinite Horizon Solution}
In order to determine the infinite horizon solution, we first notice that for $r=V$, $p^V_t =1$, $q^V_t=0$ and that the recursions \eqref{eq:recursion_and_gain} for $r = V$ are then simply given by
\begin{equation}
\begin{array}{rcl}
X^V_t &=& Q + A^\top X^V_{t+1} A + A^\top X^V_{t+1}BK^V_t \\
K^V_t &:=& \left(R+B^\top X^V_{t+1}B\right)^{-1}B^\top A,
\end{array}
\label{eq:inf_dare}
\end{equation}
that is to say the standard discrete algebraic Riccati recursion/gain.  By assumption, we have that $(X^V_t,K^V_t)\rightarrow(X^V, K^V)$, where $X^V$ and $K_V$ are, respectively, the stabilizing solution the discrete algebraic riccati equation, and the centralized LQR gain.

Now assume that $X^s_t$ is defined, and let $r\neq s \in\V$ be the unique node such that $(r,s)\in\E$.   Much as in the finite horizon case, define the following matrices:
\begin{equation}
\begin{array}{ll}
\Lambda^r = Q^{rr} +  p^r (A^{Vr})^\top X^VA^{Vr} + q^r (A^{sr})^\top X^sA^{sr} \hfill{}\\
\Psi^r = R^{rr} + p^r (B^{Vr})^\top X^VB^{Vr} + q^r (B^{sr})^\top X^sB^{sr} \hfill{} \\
\Omega^r = p^r (A^{Vr})^\top X^VB^{Vr} + q^r (A^{sr})^\top X^sB^{sr} \hfill{}
\end{array}
\end{equation}
where we have let
\begin{equation}
\begin{array}{rcl}
(p^r,q^r) = \lim_{T\to\infty}\frac{1}{T}\sum_{t=0}^{T-1}(p^r_t,q^r_t).
\end{array}
\end{equation}
Note that these limits are well defined by the assumption (\ref{eq:stationarity}).

We then have that
\begin{equation}
\begin{array}{rcl}
X^r &=& \Lambda^r + \Omega^rK^r \\
K^r &:=& - \left(\Psi^r\right)^{-1}\left(\Omega^r\right)^\top.
\end{array}
\end{equation} 

What remains to be computed is the infinite horizon average cost, which is given by (ignoring without loss the cost incurred by the uncertainty in the initial conditions)
\begin{multline}
\lim_{N\rightarrow\infty} \frac{1}{N} \sum_{t=1}^N \sum_{i=1}^2\text{Trace}\left(X^{\{i\}}_tW^i\right)\\=\sum_{i=1}^2\text{Trace}\left(X^{\{i\}}W^i\right)
\end{multline}

\section{Conclusion}
\label{sec:conclusion}
This paper presented extensions of a Riccati-based solution to a distributed
control problem with communication delays -- in particular, we now allow the communication delays to vary, but impose that they preserve partial nestedness. It
was seen that the varying delay pattern induces piecewise
linear dynamics in the state of the resulting optimal controller, with changes in dynamics
 dictated by the current \emph{effective} delay regime.

Future work will be to extend the results to systems with several players and more general delay patterns, and to remove the assumption of strong connectedness, much as was done in \cite{LL12} for the case of constant delays.  We will also seek to identify conditions on the delay process $d_t$ such that assumption (\ref{eq:stationarity}) holds.  Additionally, we will explore the setting in which the global delay regime is not known.

\bibliography{/Users/nmatni/Documents/Publications/biblio/comms,/Users/nmatni/Documents/Publications/biblio/decentralized,/Users/nmatni/Documents/Publications/biblio/matni}

\begin{appendix}
\section{Proofs}
\textbf{Proof of Lemma \ref{lem:alt_info}:}
The first two terms of \eqref{eq:alt_info} follow directly from \eqref{eq:info_recursion}. The $x^j$ component of $\I^i_t$ is then given by
\begin{multline}
\cup_{\tau=0}^t \{x^j_k\,:\, 0 \leq k \leq \tau - d^j_\tau\} = \\
\cup_{\tau=0}^t \{x^j_k\,:\, 0 \leq k \leq t - (d^j_\tau + (t-\tau))\} = \\
\{x^j_k\,:\, 0 \leq k \leq t - \min_{\tau=0,\dots,t}(d^j_\tau + (t-\tau))\} = \\
\{x_k^j\,:\, 0 \leq k \leq t-\e{t}{j}\}
\end{multline}
where the last equality follows from $d^j_t \leq D  \ \forall t\geq0$ and the definition of $\e{t}{j}$.  Noting that this is precisely the local information available to plant $j$ at time $t-\e{t}{j}$, and that the $x^i$ component of $\I^j_{t-\e{t}{j}}$ is contained in $\I^i_{t-1}\cup\{x^i_t\}$, the claim follows.
\hfill{}\QED


\textbf{Proof of Lemma \ref{lem:partial}:}
Note that $\I^i_t\subset \I^i_{t+1}$, and that $\I^i_t\subset \I^j_{t+D}$:
\begin{multline}
\I^i_t = \{x^i_{1:t}\}\cup\{x^j\,:\, 1\leq k \leq t-\e{t}{j}\} \\
\subset \{x^i_{1:t}\}\cup\{x^j\,:\, 1\leq k \leq t+D\} \\
\subset \{x^i_k\,:\, 1 \leq k \leq t+D-\e{t+D}{i}\}\cup\{x^j_{1:t+D}\} \\
= \I^j_{t+D-1}\cup\{x^j_{t+D}\}\cup \I^i_{t+D-\e{t+D}{i}} =  \I^j_{t+D}
\end{multline}
where the final inclusion follows from $\e{\tau}{i}\leq D$ for all $\tau\geq0$, and the final equalities from Lemma \ref{lem:alt_info}.  
Partial nestedness then follows from the fact that $u^i_\tau$ only affects $\I^j_t$ for $t\geq\tau+D$ due to the propagation delay between plants.  By Lemma \ref{lem:linctrl}, $u^i_t$ is a linear function of $\I^i_t$ and the same is trivially true for $x^i_t \in \I^i_t$.  We prove the final claim of the lemma by induction.

We first note that that $x_0,u_0 \in \spa{x_0} = \spa{\Ht_0}$.  We now proceed by induction, and assume that for some $t\geq0$ we have that $x_t,u_t \in \spa{\Ht_t}$.  We then have that
\begin{multline}
x_{t+1} \in \spa{\Ht_t\cup\{w_t\}}=\spa{\Ht_{t+1}} \\
u_{t+1} \in \spa{\I^1_{t+1}\cup \I^2_{t+1}} = \spa{\{x_{t+1}\}\cup\Ht_t} \\= \spa{\Ht_{t+1}}
\end{multline}\hfill{}\QED

\textbf{Proof of Lemma \ref{lem:partition}:}
(i) We begin by showing that the union in the RHS of (\ref{eq:part}) is disjoint.  This easily verified to hold for $t=0$, as all labels are the empty set except for $\La^i_0=\set{x^i_0}$.  We now proceed by induction, and suppose that the union in (\ref{eq:part}) is a disjoint one for some $t\geq 0$. We then have that
\begin{multline}
\La^V_{t+1}\cup_i\cup_{s\in\ve{t+1}{i}{-}}\La^s_{t+1} \\
= \La^V_{t}\cup_i \cup_{s\in\ve{t+1}{i}{+}}\La^s_t\cup_{s\in\ve{t+1}{i}{--}}\La^s_{t}\cup\La^i_{t+1}
\label{eq:disj}
\end{multline}
where the equality follows from simply applying the recursion rules \eqref{eq:labelrecursions} and Lemma \ref{lem:shrink}.  We first note that by the induction hypothesis, $\La^V_{t}\cap\cup_i\cup_{s\in\ve{t+1}{i}{--}}\La^s_{t} = \emptyset$.  Additionally, by construction, we have that $\cup_i \cup_{s\in\ve{t+1}{i}{+}}\La^s_t\cap\cup_i\cup_{s\in\ve{t+1}{i}{--}}\La^s_{t}=\emptyset$.  We note that $\La^i_{t+1} = \set{w^i_t}$ is the new information available at time $t+1$, and thus $\La^i_{t+1}\cap\La^s_t=\emptyset$ for all $s\in\V$.  Finally, noting that for all $\La^1_{t+1}\cap\La^2_{t+1}=\emptyset$, we have that \eqref{eq:disj} is a disjoint union, proving the claim.

It now suffices to show that (\ref{eq:part}) is also a covering of the noise history.  To that end, notice that for $t=0$, this follows immediately from $\La^i_0 = \{x^i_0\}$, and $\Ht_0 = \{x_0\}$.  Now suppose that (\ref{eq:part}) is a covering for some $t\geq 0$.  We then have that
\begin{multline}
\Ht_{t+1} = \Ht_t\cup_i \La^i_{t+1} = \La^V_t\cup_i\cup_{s\in\ve{t}{i}{-}}\La^s_t\cup\La^i_{t+1} \\
= \La^V_t\cup_i\cup_{s\ni i,\, |s|\leq \e{t}{i}+1 }\La^s_{t+1} \hfill{} \\
= \La^V_t\cup_i\cup_{s\in\ve{t+1}{i}{-}}\La^s_{t+1}\cup_{s' \ni i,\, \e{t+1}{i}<|s'|\leq\e{t}{i}+1}\La^{s'}_{t+1} \\
= \La^V_{t+1}\cup_i\cup_{s \in \ve{t+1}{i}{-}}\La^s_{t+1}. \hfill{}
\end{multline}
The third equality follows from applying the induction hypothesis, the fourth by applying the recursion rules for the label sets, and the before last equality from noticing that $\e{t+1}{i}\leq \e{t}{i}+1$.  To prove the final equality, it suffices to show that $\La^V_t \cup_i \cup_{s' \ni i,\, \e{t+1}{i}<|s'|\leq\e{t}{i}+1}\La^{s'}_{t+1} = \La^V_{t+1}$.  This follows by applying the recursion rules and Lemma \ref{lem:shrink} as follows:
\begin{multline}
\La^V_t \cup_i \cup_{s' \ni i,\, \e{t+1}{i}<|s'|\leq\e{t}{i}+1}\La^{s'}_{t+1} \\
= \La^V_t \cup_i \cup_{s' \ni i} \cup_{\e{t+1}{i}\leq|s'|\leq\e{t}{i}}\La^{s'}_{t}\cup_{|s'|\geq \e{t}{i}}\La^{s'}_{t-1} \\
= \La^V_t \cup_i \cup_{s' \ni i} \cup_{\e{t+1}{i}\leq|s'|\leq\e{t}{i}}\La^{s'}_{t}\cup_{|s'|\geq \e{t}{i}+1}\La^{s'}_{t} \\
= \La^V_t \cup_i \cup_{s \in \ve{t+1}{i}{+}}\La^s_t = \La^V_{t+1}
\end{multline}

(ii) We proceed by induction once again.  This holds trivially for $t=0$.  Now suppose it to be true for some $t\geq 0$.  We have that $\I^i_{t+1} = \I^i_t\cup\I^j_{t-(\e{t+1}{j}-1)}\cup\{x^i_{t+1}\}$.  Taking the linear span of both sides, we then obtain 
\begin{multline}
 \spa{\I^i_{t+1}} = \spa{\I^i_t}+\spa{\I^j_{t-(\e{t+1}{j}-1)}} + \spa{w^i_t} \hfill{}\\
= \spa{\La^V_t} + \sum_{s\in\ve{t}{i}{-}}\spa{\La^s_t} + \dots \hfill{} \\ \sum_{r\in\ve{t+1}{j}{--}}\spa{\La^r_{t-(\e{t+1}{j}-1)}} +  \spa{\La^i_{t+1}} 
\label{eq:partii_interim}
\end{multline}

By the same arguments used in the second part of the proof of part (i), we have that $\spa{\sum_{s\in\ve{t}{i}{++}}\La^s_t}\subset \spa{\La^V_t}$. Also notice that applying the recursion for $\La^s_{t+1}$ to the $\La^r_{t-(\e{t+1}{j}-1)}$ term $\e{t+1}{j}-1$ times, and that for $r\rightarrow \dots \rightarrow s'$, we have that $|s'| = |r| + \e{t+1}{j}-1\geq \e{t+1}{j}$.  We may  then write \eqref{eq:partii_interim} as
\begin{multline}
\spa{\La^V_t} + \sum_{s\ni i} \spa{\La^s_t} + \sum_{s'\in\ve{t+1}{j}{+}}\spa{\La^{s'}_{t}} \\ 
= \spa{\La^V_t} + \sum_{k=1}^2\sum_{s\in\ve{t+1}{k}{+}}\spa{\La^s_t} + \dots \\ \sum_{s \in \ve{t+1}{i}{--}}\spa{\La^s_t} + \spa{\La^i_{t+1}}.
\label{eq:partii_almost}
\end{multline}
The first two terms of the final equality are precisely the expression for $\spa{\La^V_{t+1}}$, whereas the final two terms may be combined by applying the recursion rules to the summation, yielding $\sum_{s\in\ve{t+1}{i}{-}}\spa{\La^s_{t+1}}$.   We therefore have that \eqref{eq:partii_almost} is equal to
\begin{multline}
\spa{\La^V_{t+1}} + \sum_{s\in\ve{t+1}{i}{-}}\spa{\La^s_{t+1}} =\\ \spa{\La^V_{t+1} \cup_{s\in \ve{t+1}{i}{-}}\La^s_{t+1}}
\end{multline}
proving the claim.
\hfill{}\QED

\textbf{Proof of Lemma \ref{lem:dynamics}:}
The recursive nature of the label sets ensure that $\zeta^s \in \spa{\La^s_t}$ for all $t\geq0$.  Thus it suffices to show that these dynamics preserve the state decomposition (\ref{eq:decomp}).
\begin{multline}
\zeta^V_{t+1} + \sum_{i=1}^2\sum_{s\in\ve{t+1}{i}{-}} I^{V,s}\zeta^s_{t+1} \hfill{}\\ 
= A\zeta^V_t + B\varphi^V_t +  \sum_{i=1}^2 \sum_{r\in\ve{t+1}{i}{++}}(A^{Vr}\zeta^r_t + B^{Vr}\varphi^r_t) \\
+\sum_{i=1}^2\sum_{s\in\ve{t+1}{i}{-}}I^{V,s}\left( A^{sr}\zeta^r_t + B^{sr}\varphi^r_t\right) +w_t\\
= A\left(\zeta^V_{t} + \sum_{s\in\V}I^{V,s}\zeta^s_t\right)+B\left(\varphi^V_{t} + \sum_{s\in\V}I^{V,s}\varphi^s_t\right) +w_t\\
= A\left(\zeta^V_t + \sum_{i=1}^2\sum_{s\in\ve{t}{i}{-}} I^{V,s}\zeta^s_t\right) \hfill{}\\+B\left(\varphi^V_t + \sum_{i=1}^2\sum_{s\in\ve{t}{i}{-}} I^{V,s}\varphi^s_t\right) +w_t\\
= Ax_t + Bu_t + w_t = x_{t+1}
\end{multline}
where the first equality followed from applying the update dynamics (\ref{eq:dynamics}), and the third from noting that certain components of the state and control decomposition are zero due to the effective delays seen by the controllers.  The fourth equality follows from equation (\ref{eq:decomp}), and the final one from (\ref{eq:statespace}).
\hfill{}\QED

\end{appendix}
\end{document}